\documentclass[12pt,reqno]{amsart}

\usepackage{setspace} \usepackage[dvips]{graphicx, hyperref}

\newcommand\N {{\mathbb N}}

\author[T.B.]{Tiziana Bascelli} \address{T. Bascelli, Via
S. Caterina, 16, 36030 Montecchio P.no (VI), Italy}
\email{tiziana.bascelli@virgilio.it}

\author[P.B.]{Piotr B\l{}aszczyk}\address{P. B\l{}aszczyk, Institute
of Mathematics, Pedagogical University of Cracow,
Poland}\email{pb@up.krakow.pl}

\author[V.K.]{Vladimir Kanovei} \address{V. Kanovei, IPPI, Moscow,
and MIIT, Moscow, Russia}\email{kanovei@googlemail.com}

\author[K.K.]{Karin U. Katz}\address{K. Katz, Department of
Mathematics, Bar Ilan University, Ramat Gan 52900 Israel}
\email{katzmik@math.biu.ac.il}

\author[M.K.]{Mikhail G. Katz}\address{M. Katz, Department of
Mathematics, Bar Ilan University, Ramat Gan 52900 Israel}
\email{katzmik@macs.biu.ac.il}

\author[D.Sc.]{David M. Schaps}\address{D. Schaps, Department of
Classical Studies, Bar Ilan University, Ramat Gan 52900 Israel}
\email{dschaps@mail.biu.ac.il}

\author[D.Sh.]{David Sherry}\address{D. Sherry, Department of
Philosophy, Northern Arizona University, Flagstaff, AZ 86011, US}
\email{David.Sherry@nau.edu}

\begin{document}

%\doublespacing

\thispagestyle{empty}

\title[Leibniz vs Ishiguro] {Leibniz vs Ishiguro: Closing a
quarter-century of syncategoremania}

\begin{abstract}
Did Leibniz exploit infinitesimals and infinities \emph{\`a la
rigueur}, or only as shorthand for quantified propositions that refer
to ordinary Archimedean magnitudes?  Chapter~5 in (Ishiguro 1990) is a
defense of the latter position, which she reformulates in terms of
Russellian \emph{logical fictions}.  Ishiguro does not explain how to
reconcile this interpretation with Leibniz's repeated assertions that
infinitesimals violate the Archimedean property, viz., Euclid's
\emph{Elements}, V.4.  We present textual evidence from Leibniz, as
well as historical evidence from the early decades of the calculus, to
undermine Ishiguro's interpretation.  Leibniz frequently writes that
his infinitesimals are \emph{useful fictions}, and we agree; but we
shall show that it is best not to understand them as logical fictions;
instead, they are better understood as \emph{pure fictions}.

\end{abstract}

\maketitle
\tableofcontents

\section{Logical fictions}
\label{one}
\label{honor1}

If a publisher were to announce to the public that in addition to its
fiction titles, it offers a variety of cookbooks as well, no one would
interpret this as meaning that the fiction titles turn out to be
cookbooks in disguise when their content is properly clarified and
made explicit.  Yet when Leibniz announced that ``il ne faut pas
s'imaginer que la science de l'infini est \ldots{} reduite \`a des
fictions; car il reste tousjours un infini syncategorematique''%
\footnote{This passage is discussed in more detail in
Subsection~\ref{71}.}
\cite[p.~93]{Le02a}, H.~Ishiguro proposed just this type of
interpretation of Leibnizian fictions in terms of a Weierstrassian
cookbook.  
%
%To put it bluntly, while Leibniz talks about both apples and oranges,
%Ishiguro claims that Leibnizian apples turn out to be oranges when
%properly peeled.

Twenty-five years ago Ishiguro presented her interpretation of
Leibnizian infinitesimals as logical fictions \cite[Chapter~5]{Is90}.
Ishiguro's interpretive strategy employs what Russell called logical
or symbolic fictions \cite [pp.~45 and 184]{Ru19}.

Her analysis has not been seriously challenged.  The situation has
reached a point where the literature contains statements of Ishiguro's
hypothesis as fact, without any further attribution, as in the
following:
\begin{quote}
Robinson's infinitesimal is a static quantity, whereas Leibniz's
infinitesimals are ``syncategorematic,'' i.e., they are as small as is
necessary, such that there is always a quantity that is smaller than
the smallest given quantity. \cite[p.~15]{Du13}
\end{quote}
Such a syncategorematic reading of Leibnizian infinitesimals is
endorsed by Leibniz scholars Arthur, Goldenbaum, Knobloch, Levey,
Nachtomy, and others, as detailed in Subsection~\ref{22} below.
Recent advances in Leibniz scholarship suggest the time has come to
re-evaluate her interpretation.  This text presents a number of
difficulties for the thesis that Leibnizian infinitesimals are logical
fictions.

\subsection{Defending Leibniz's honor}
\label{11}

The context for Ishiguro's analysis was the general sense that no
appeal to infinitesimals \emph{\`a la rigueur} could stand
philosophical scrutiny.  More specifically, her reading is based on
the premise that prior interpretations of Leibnizian infinitesimals,
in the spirit of the infinitesimals of Bernoulli and Euler, must
surely involve confusion or even logical inconsistency.  This premise
is spelled out in the title of her text ``La notion dite confuse de
l'infinit\'esimal chez Leibniz" \cite{Is86}, an early version of her
Chapter~5.  As she writes there, ``This is because the concept of
infinitesimal was seen as being \emph{confused}.''  \cite[p.~79]{Is90}
[emphasis added] Furthermore,
\begin{quote}
The second kind of critic acknowledges that Leibniz was interested in
foundational issues, but after examination sees basic
\emph{inconsistencies} in his views.  (ibid., p.~80) [emphasis added]
\end{quote}
Thus, Ishiguro purports to defend Leibniz's honor as an
\emph{unconfused} and \emph{consistent} logician by means of her
syncategorematic reading.  Meanwhile, in the first edition of her
book, Ishiguro wrote:
\begin{quote}
Leibniz's philosophy of logic and language makes far more sense in
every aspect than has generally been thought, let alone that his
thought is more coherent than Russell allowed. \cite[p.~16]{Is72}
\end{quote}
We argue that such an appreciation of Leibniz applies equally well to
his infinitesimal calculus in the spirit of Bernoulli and Euler.
Ishiguro goes on to write:
\begin{quote}
In many respects, it is much less dated than the theories of Locke and
Berkeley, and even of Kant.  (ibid.)
\end{quote}
According to Ishiguro, the superiority of Leibniz's thought over that
of Locke and Berkeley manifests itself also in Leibniz's rejection of
empiricism.%
\footnote{\label{f2}Thus, Leibniz's disagreements with empiricism are
mentioned in the final paragraph of Ishiguro's 1972 book:
``[Leibniz's] disagreement with many of the views of the empiricists,
as with those of the Cartesians, sprang from his belief that their
theories failed to account for the complex facts which fascinated him,
whether these were about the language we have or about the concepts we
use.''  (ibid., p.~145) These comments on empiricism and Cartesianism
appeared at the end of section~6 entitled ``Concepts resolvable at
infinity" in the final chapter~7, entitled ``Necessity and
Contingency.''  Ishiguro's section~6 ``Concepts resolvable at
infinity" is still present in the second (1990) edition, though
``Necessity and contingency" is now chapter~9 rather than 7 (this is
due in part to the addition of Chapter~5 seeking to reduce
infinitesimals to quantified propositions).  The comment on empiricism
and Cartesianism disappeared from the second edition, but here
Ishiguro writes that Leibniz's reasoning ``is not of an empiricist
kind like that of Berkeley.''  \cite[p.~85]{Is90}}

We similarly believe that Leibniz was not confused and likewise intend
to defend his honor, in this case \emph{against} Ishiguro's reading.
%
%{honor1}
We will see that at a few key junctures, Ishiguro is forced to defend
her reading by attributing confusion to Leibniz; see
Subsections~\ref{honor2} and~\ref{honor4}.  On at least one occasion,
Ishiguro misrepresents what Leibniz wrote so as to buttress her
position; see Subsection~\ref{honor3}.
We shall argue that the appeal to \emph{logical} fictions is neither
necessary to defend Leibniz's honor, nor warranted in view of the
actual content of Leibniz's mathematics and philosophy.

\subsection{Categorematic vs syncategorematic} 
\label{12}

According to Ishiguro, expressions like~$dy/dx$, which appear to refer
to infinitesimals, are not in fact referring, denoting, or
\emph{categorematic}, expressions.  Rather, they are
\emph{syncategorematic} expressions, namely expressions which
disappear when the logical content of the propositions in which they
occur is properly clarified and made explicit.  Writes Ishiguro:
\begin{quote}
The word `infinitesimal' does not designate a special kind of
magnitude.  In fact, it does not designate%
\footnote{Ishiguro uses \emph{designate} as an intransitive verb, and
similarly for the verbs \emph{denote} and \emph{refer}.  A term is
said not to \emph{refer} when the term does not actually refer to
anything but rather is awaiting a clarification of the logical content
of the sentence it occurs in, which would make the term disappear.  An
example is provided by Weierstrass's use of the term
\emph{infinitesimal} as discussed in Subsection~\ref{21}.}
at all. \cite[p.~83]{Is90}
\end{quote}
A few pages later, she clarifies the nature of her non-designating
claim in the following terms:
\begin{quote}
we can paraphrase the proposition with a universal proposition with an
embedded existential claim. (ibid., p.~87)
\end{quote}
In conclusion,
\begin{quote}
Fictions [such as Leibnizian infinitesimals] are \emph{not entities}
to which we refer.  \ldots{} They are correlates of ways of speaking
which can be reduced to talk about more standard kinds of entities.
(ibid., p.~100) [emphasis added]
\end{quote}
Such fictions (which are \emph{not entities}) are exemplified by
Leibnizian infinitesimals, in Ishiguro's view.  Her contention is
that, when Leibniz talked about infinitesimals, what he \emph{really
meant} was a certain quantified proposition, or more precisely a
quantifier-equipped proposition.  In short, Leibniz was talking about
ordinary numbers.

For the 17th century context see \cite{Al14}.  Ishiguro does mention
``Leibniz's followers like Johann Bernoulli, de l'Hospital, or Euler,
who were all brilliant mathematicians rather than philosophers,''
\cite[pp.~79-80]{Is90} but then goes on to yank Leibniz right out of
his historical context by claiming that their \emph{modus operandi}
\begin{quote}
is prima facie a strange thing to ascribe to someone who, like
Leibniz, was obsessed with general methodological issues, and with the
logical analysis of all statements and the well-foundedness of all
explanations.  (ibid., p.~80)
\end{quote}
Having thus abstracted Leibniz from his late 17th century context,
Ishiguro proceeds to insert him in a late 19th century Weierstrassian
cookbook.  Such an approach to a historical figure would apparently
not escape Unguru's censure:
\begin{quote}
It is \ldots{} a historically unforgiveable sin \ldots{} to assume
wrongly that mathematical equivalence is tantamount to historical
equivalence.  \cite[p.~783]{Un76}
\end{quote}
Ishiguro seems to have been aware of the problem and at the end of her
Chapter~5 she tries again to explain
\begin{quote}
why I believe that Leibniz's views on the contextual definition of
infinitesimals is [sic] different from those of other mathematicians
of his own time who sought for operationist definitions for certain
mathematical notions \cite[p.~99]{Is90},
\end{quote}
but with limited success.

\section{Testing the limits of syncategorematics}
\label{ch}

We take it that `infinitesimal' expressions do designate insofar as
our symbolism allows us to think about infinitesimals.  It should be
emphasized that our contention that a Leibnizian infinitesimal
\emph{does} designate does not imply that they designate entities on a
par with monads, material objects, or ideal entities.  While
\emph{infinitesimal} is a designating expression for Leibniz, it
designates a fictional entity.  Likewise, for Leibniz, \emph{imaginary
quantity} designates a fictional entity (see Subsection~\ref{7june}).
The literature contains a considerable amount of confusion on this
subject, as in the following:
\begin{quote}
The use of fictitious quantities could lead to the erroneous idea of
objects whose existence is assured by their very definition and
therefore to ascribing a modern conception to Leibniz. In reality,
what finds its foundation in Nature cannot be created by the human
mind by means of a definition. \cite[p.~36]{Fe08}
\end{quote}
Now the matter of creating by definition is a tricky one.  Leibniz
certainly denies that definitions carry existential commitments.  In
Leibniz, infinitesimals are created at the syntactic level by
postulation, which has a subtle relation to existence.  This must be
so, since the difficulty Ferraro perceives arises equally for real
numbers.  The article \cite[p.~322]{Le95b} introduces infinitesimals
specifically by invoking a definition, namely Euclid Def.~V.4, and
postulating that infinitesimals are entities that fail to satisfy the
latter.

In order to test the range of applicability of Ishiguro's
syncategorematic reading, we consider the following two extreme cases.

\subsection{Weierstrass on infinitesimals}
\label{21}

On the one hand, there does exist a context where Ishiguro's
\emph{logical fiction} hypothesis may be on solid ground.  On
occasion, Weierstrass mentions an infinitesimal definition of
continuity.  This is Cauchy's original definition of continuity (see
\cite[p.~34]{Ca21}) of a function~$y=f(x)$:
\emph{infinitesimal~$x$-increment always produces an infinitesimal
change in~$y$}.  Thus, Weierstrass wrote:
\begin{quote}
Finally, once the concept of the infinitely small has been grasped
correctly, one can define the concept of the continuity of a function
in the vicinity of $a$ as follows: that infinitely small changes in
the arguments correspond to infinitely small changes in the value of
the function in the vicinity of $a$.%
\footnote{In the original: ``Endlich kann man, den Begriff des
unendlich Kleinen richtig gefa\ss{}t, den Begriff der Stetigkeit einer
Funktion in der N\"ahe von a auch dadurch definieren, da\ss{}
unendlich kleinen \"Anderungen der Argumente unendlich kleine
\"Anderungen des Funktionswertes in der N\"ahe von a entsprechen
sollen.''}
\cite[p.~74]{We86}
\end{quote}
It may be reasonable to conjecture that when Weierstrass refers to an
infinitesimal, he \emph{always} means (unlike Leibniz, on our reading)
a kind of logical fiction.  Here an infinitesimal is shorthand for a
longer paraphrase expressed by a proposition whose quantifiers range
over ordinary real numbers, namely the sort of proposition that
typifies Weierstrass's contribution to the foundations of analysis.

On the other hand (and at the other extreme), an infinitesimal is
\emph{not} meant to be a shorthand for a quantified paraphrase in the
context of modern infinitesimal frameworks such as those of Robinson
(see \cite{Ro61}), Lawvere-Kock (see \cite{Ko06}), or J. Bell
\cite{Be06}.
%
%\footnote{However, see footnote~\ref{nelson} for a discussion of
%conservativity.}
%
Note that Robinson as a formalist distanced himself from Platonist and
foundationalist views: ``mathematical theories which, allegedly, deal
with infinite totalities do not have any detailed \ldots{}
reference.''  \cite[p.~42]{Ro75}

Ishiguro's argument is based on first philosophical principles (rather
than on historical analysis or careful textual study) that are so
general that, while it might apply to Weierstrass, it is difficult to
see what would prevent her from applying it to Robinson, as well.  Yet
scholars agree that Robinson's infinitesimals are not logical
fictions; nor is his continuum Archimedean.

\subsection{Syncategorematic \emph{vs} fictionalist}
\label{22}

The syncategorematic interpretation of Leibnizian infinitesimals is
the starting point of much recent Leibniz scholarship.  Leading
Leibniz scholar E. Knobloch writes:
\begin{quote}
To my knowledge most of the historians of mathematics are convinced
that Leibniz used an Archimedean continuum: Leibniz himself referred
to the Greek authority in order to justify his procedure. \cite{Kn14}
\end{quote}
Writes U.~Goldenbaum:
\begin{quote}
That Leibniz as a mature mathematician and philosopher did not take
infinitesimals to be real entities, but rather as finite quantities,
was clarified as early as 1972 by Hid\'e Ishiguru [sic]. \cite[p.~76,
note~59]{Go08}.
\end{quote}
\cite[note~25]{Ra15} similarly endorses Ishiguro's Chapter~5.  Both
S.~Levey and \cite[p.~20]{Ar08}, \cite[p.~554]{Ar13} take Ishiguro's
interpretation as settled and have concentrated, instead, on
demonstrating that Leibniz embraced the syncategorematic
interpretation of infinitesimals as early as 1676.  The following,
from Levey, is typical:

\begin{quote}
[B]y April of 1676, with his early masterwork on the calculus,
\emph{De Quadratura Arithmetica}, nearly complete, Leibniz has
\emph{abandoned an ontology of actual infinitesimals} and adopted the
syncategorematic view of both the infinite and the infinitely small as
a philosophy of mathematics and, correspondingly, he has arrived at
the official view of infinitesimals as fictions in his
calculus. \cite[p.~133]{Le08} [emphasis added]
\end{quote}
O.~Nachtomy chimes in:
\begin{quote}
\ldots{} Richard Arthur makes a very convincing argument that
Leibniz's syncategorematic view of infinitesimals was developed in the
very early 1670s and matured in 1676. \cite{Na09}
\end{quote}

We don't intend to disagree with Levey and others that Leibniz may
have ``abandoned an ontology of actual infinitesimals'' early on.
However, we object to the conflation of the views of the
syncategorematicist and the fictionalist.  The syncategorematic
interpretation is a fictionalist interpretation, to be sure, but the
converse is not the case.

In what follows we shall demonstrate that Leibniz understood this
point and had good reason to embrace a different variety of
fictionalism, which we call \emph{pure} fictionalism.%
\footnote{Perhaps we may be allowed to quote Leibniz's own description
of his method: \emph{My arithmetic of infinites is pure, Wallis' is
figurate} (``Arithmetica infinitorum mea est pura, Wallisii figurata''
\cite[p.~102]{Le72}), as translated in \cite[p.~46]{Be08}.  Note that
the text in question, \emph{De progressionibus et de arithmetica
infinitorum}, predates the \emph{Arithmetic Quadrature}.  What Leibniz
meant by \emph{figurata} is not entirely clear, nor does the context
offer an indication.  The most likely interpretation seems to be that
Wallis relied on induction from geometric \emph{figures}, which
Leibniz' method did not require.}
Modern exponents of this variety of fictionalism include Hilbert and
Robinson (see \cite{KS1}).  

That Leibniz considered an alternative version of fictionalism will
come as a surprise mainly to scholars whose outlook presumes that the
epsilon-delta style of analysis, promoted by the ``triumvirate'' of
Cantor, Dede\-kind, and Weierstrass (see \cite[p.~298]{Bo49}) is the
embodiment of inevitable progress climaxing in the establishment of
the foundations of real analysis purged of infinitesimals; related
issues are explored in articles \cite{KK15}, \cite{Ka15a},
\cite{Ka15b}.

\subsection{Summary of Ishiguro's hypothesis}
\label{23}

According to Ishiguro, Leibniz's conception of continuity (i.e., the
continuum) is Archimedean.  On the syncategorematic reading, talk
about infinitesimals involves only expressions which do not denote
anything.  The position as expressed in \cite[Chapter~5]{Is90} can
therefore be summarized in terms of the following three contentions.
\begin{enumerate}
\item
Taking Leibnizian infinitesimals at face value requires one to see
Leibniz as confused (ibid., p.~79) and/or inconsistent (ibid., p.~80);
\item
A term that seems to express a Leibnizian infinitesimal does not
actually designate, denote, or refer, and is a logical fiction;
\item
Leibniz's continuum is Archimedean.
\end{enumerate}

None of these can be sustained in light of Leibniz's philosophical and
mathematical texts.

\section{Analysis of Ishiguro's contentions}
\label{three}

Let us analyze Ishiguro's hypothesis as summarized in
Subsection~\ref{23}.  Ishiguro's contention (1) concerns scholars like
C.~Boyer and J.~Earman (see \cite[p.~80]{Is90}).  However, a perusal
of their work reveals that the ultimate source of the inconsistency
claim is Berkeley's \emph{departed quantities}.  Thus, Ishiguro's
contention~(1) echoes Berkeley's claim that inconsistent properties
have been attributed to~$dx$, viz.,
\[
(dx\not=0)\wedge (dx=0).%
\footnote{Berkeley's critique was dissected into its logical and
metaphysical components in \cite{Sh87}.  The logical criticism
concerns the alleged inconsistency expressed by the conjunction
$(dx\not=0)\wedge(dx=0)$, while the metaphysical criticism is fueled
by Berkeley's empiricist doubts about entities that are below any
finite perceptual threshold; see Subsection~\ref{11} and
note~\ref{f2}.}
\]
However, Berkeley's claim ignores Leibniz's generalized relation of
equality (see Subsection~\ref{32}).

\subsection{Of $dx$ and $(d)x$}
\label{31}

Ishiguro does not appear to be an attentive reader of H.~Bos, and in
fact on page 81 she misrepresents his position.  She quotes Bos to the
effect that Leibniz eventually introduced the finite (assignable)
differentials.  For these he used the notation~$(d)x$ in place
of~$dx$.  These satisfy the equality on the nose
\[
(d)y=L\;(d)x,
\]
where $L$ is what we would call today the derivative.  Bos does say
that.

However, Ishiguro further implies that, according to Bos, these
$(d)x$'s completely replaced the earlier~$dx$'s.  She then goes on to
disagree with her strawman Bos by claiming that Leibniz never changed
his mind about infinitesimals (namely that they were always logical
fictions).  On p.~81 Ishiguro writes: ``Bos talks (perhaps naturally
as a post-Robinsonian) as if it is quite clear what it means for a
magnitude to be infinitely small, and that Leibniz first assumed the
existence of such things.'' Bos may well be surprised to find himself
described as a post-Robinsonian, especially given what he wrote about
Robinson in \cite[Appendix 2]{Bos}; see further in
Subsection~\ref{33}.

Contrary to Ishiguro, Bos never asserted that the~$dx$'s disappeared
with the introduction of~$(d)x$'s.  Bos merely reports that Leibniz
introduced the additional concepts~$(d)x$, not that they completely
replaced the~$dx$'s, which they certainly never did.  Thus, the late
piece \emph{Cum Prodiisset} \cite{Le01c} features both the $(d)x$'s
and the~$dx$'s, as well as the crucial distinction between
\emph{assignable} and \emph{inassignable}:
\begin{quote}
\ldots although we may be content with the assignable quantities
$(d)y$, $(d)v$, $(d)z$, $(d)x$, etc., \ldots{} yet it is plain from
what I have said that, at least in our minds, the unassignables
[\emph{inassignables} in the original Latin] $dx$ and $dy$ may be
substituted for them by a method of supposition even in the case when
they are evanescent; \ldots{} (as translated in \cite[p.~153]{Ch}).
\end{quote}

\subsection{Law of homogeneity}
\label{32}

Bos notes that Leibniz already mentioned his \emph{law of homogeneity}
in \emph{Nova Methodus} \cite{Le84}.  Leibniz explained the law in a
letter to Wallis \cite[p.~63]{Le99}, and gave the most detailed
presentation in his 1710 piece \cite{Le10b} mentioning the
\emph{transcendental law of homogeneity} (TLH) in the title.  The law
involves, roughly, discarding higher-order terms.

Leibniz was using the relation of equality in a generalized sense of
\emph{equality up to}, as mentioned in his \emph{Responsio} (see e.g.,
Subsection~\ref{honor3} below, sentences labeled [1] and [2]).  This
means that~$dx$ does not turn out to be zero at the end of the
calculation but, rather, is discarded in an application of TLH.
\emph{Equality up to} undermines the claim of logical inconsistency
(alleged by Berkeley) without a need to dip into a Weierstrassian
cookbook with hidden quantifiers \`a la Frege.

An antecedent to the Leibnizian generalized equality is found in
Fermat's relation of \emph{adequality}; see \cite{Ba13}, \cite{Ba14},
\cite{Ci}, \cite{KSS13}.  Leibniz in fact mentions Fermat's method in
the context of a discussion of the generalized notion of equality.
Here Leibniz is objecting to Nieuwentijt's postulation that the square
of an infinitesimal term should be exactly nothing:
\begin{quote}
Quod autem in aequationibus Fermatianis abjiciuntur termini, quos
ingrediuntur talia quadrata vel rectangula, non vero illi quos
ingrediuntur simplices lineae infinitesimae, ejus ratio non est, quod
hae sint aliquid, illae vero sint nibil, sed quod termini ordinarii
per se destruuntur, hinc restant tum termini, quos ingrediuntur Iineae
simplices infinite parvae, tum quos ingrediuntur harum quadrata vel
rectangula: cum vero hi termini sint illis incomparabiliter minores,
abjiciuntur.  Quod si termini ordinarii non evanuissent, etiam termini
infinitesimarum linearum non minus, quam ab his quadratorum abjici
debuissent. \cite[p.~323]{Le95b}
\end{quote}
We translate this as follows:
\begin{quote}
But the reason that in Fermat's equations, terms incorporating squares
or similar products are discarded, but not those containing simple
infinitesimal lines [i.e., segments], is not that the latter are
something, whereas the former are, on the contrary, nothing; but
rather that the ordinary terms cancel each other out, whence there
then remain terms containing infinitely small simple lines, and also
those containing their squares or products: but since the latter terms
are incomparably smaller than the former, they are discarded.  Because
if the ordinary terms did not disappear, then the terms of the
infinitesimal lines would have to be discarded no less than their
squares or products.%
\footnote{A French translation is in \cite[p.~329]{Le1989}.}
\end{quote}
Here Leibniz describes Fermat's method in a way similar to Leibniz's
own.  

%To follow Ishiguro's Russellian (but non-Unguruan) procedure to its
%logical conclusion, one would need similarly to abstract Fermat from
%his historical context and insert him in a 19th century Weierstrassian
%cookbook (see Subsection~\ref{12}).

Ishiguro's contention (2) is based on the \emph{fictional} status of
Leibnizian infinitesimals.  To be sure, Leibniz frequently describes
his infinitesimals as ``useful fictions.''  However, their fictional
nature could merely mean to Leibniz that they lack reference to either
a material object or an ideal entity, as Leibniz often writes, not
necessarily that they are \emph{logical} fictions as Ishiguro claims.

Ishiguro's contentions (2) and (3) amount to a claim of
proto-Weier\-strassian hidden quantifier ranging over ordinary
Archimedean quantities.  One of the passages claimed to support such a
reading of Leibniz is a letter to Pinson dated 29 august 1701, where
Leibniz writes:
\begin{quote}
[I]n lieu of the infinite or infinitely small, we take quantities as
great or as small as it is required so that the error would be less
than the given error such that we do not differ from the style of
Archimedes except in the expressions [....] (as translated in
\cite[p.~71]{Th12})
\end{quote}
(we have retained Tho's precise punctuation which turns out to be
significant; see below).  This passage is an optimistic expression of,
in Jesseph's phrase, a \emph{grand programmatic statement}; see
Section~\ref{four}.  We will analyze this passage further in
Section~\ref{five}.

\subsection{Ishiguro, Bos,  Robinson}
\label{33}

Given Ishiguro's \emph{post-Robinso\-nian} description of Bos (see
Subsection~\ref{31}), it will prove instructive to examine the matter
in more detail.  On the one hand, Robinson famously argued for
continuity between the Leibnizian framework and his own.  On the
other, Bos rejected such claims of continuity in his Appendix~2:
\begin{quote}
\ldots{} the most essential part of non-standard analysis, namely the
proof of the existence of the entities it deals with, was entirely
absent in the Leibnizian infinitesimal analysis, and this constitutes,
in my view, so fundamental a difference between the theories that the
Leibnizian analysis cannot be called an early form, or a precursor, of
non-standard analysis \cite[p.~83]{Bos}.
\end{quote}

Bos's comment is not sufficiently sensitive to the dichotomy of
practice (or procedures) versus ontology (or foundational account for
the entities such as numbers).  While it is true that Leibniz's
calculus contains nothing like a set-theoretical existence proof,
nonetheless there do exist Leibnizian procedures exploiting
infinitesimals that find suitable proxies in the procedures in the
hyperreal framework.  In other words, there are close formal analogies
between inference procedures in the Leibnizian calculus and the
Robinsonian calculus.  See \cite{Re13} for a related discussion in the
context of Euler.  The relevance of such hyperreal proxies is in no
way diminished by the fact that set-theoretic \emph{foundations} of
the latter (``proof of the existence of the entities,'' as Bos put it)
were obviously as unavailable in the 17th century as set-theoretic
foundations of the real numbers.

In the context of his discussion of ``present-day standards of
mathematical rigor'', Bos writes: 
\begin{quote}
it is understandable that for mathematicians who believe that these
present-day standards are final, nonstandard analysis answers
positively the question whether, after all, \emph{Leibniz was right}
\cite[p.~82, item~7.3]{Bos}. [emphasis added]
\end{quote}
The context of the discussion makes it clear that Bos's criticism
targets Robinson.  If so, Bos's criticism suffers from a strawman
fallacy, for Robinson specifically wrote that he did not consider set
theory to be the foundation of mathematics, and being a formalist, he
did not subscribe to the view attributed to him by Bos that
``present-day standards are final.''  Robinson expressed his position
on the status of set theory as follows:
\begin{quote}
an infinitary framework such as set theory \ldots{} cannot be regarded
as the ultimate foundation for mathematics \cite[p.~45]{Ro69}; see
also \cite[p.~281]{Ro66}.
\end{quote}
Furthermore, contrary to Bos's claim, Robinson's goal was not to show
that ``Leibniz was right'' as Bos claimed.  Rather, Robinson sought to
provide hyperreal proxies for the inferential procedures commonly
found in Leibniz as well as Euler and Cauchy; for the latter see e.g.,
\cite{BK}.  Leibniz's procedures, involving as they do infinitesimals
and infinite numbers, seem far less puzzling in light of their B-track
hyperreal proxies than from the viewpoint of the received A-track
frameworks; see Section~\ref{four}.

Some decades later, Bos has distanced himself from his Appendix~2 in
the following terms (in response to a question from one of the authors
of the present text):
\begin{quote}
An interesting question, what made me reject a claim some 35 years
ago?  I reread the appendix [i.e., Appendix~2] and was surprised about
the self assurance of my younger self.  I'm less definite in my
opinions today - or so I think.  \cite{Bos10}
\end{quote}
And he continues:
\begin{quote}
You're right that the appendix was not sympathetic to Robinson's view.
Am I now more sympathetic?  If you talk about ``historical continuity"
I have little problem to agree with you, given the fact that one can
interpret continuity in historical devlopments in many ways; even
revolutions can come to be seen as continuous developments. (ibid.)
\end{quote}
While Bos acknowledges that his Appendix~2 was ``unsympathetic to
Robinson's view'' we must also point out that his opinions as
expresssed in Appendix~2 were based on mathematical misunderstandings
(particularly in connection with the transfer principle, as discussed
in \cite[Section~11.3]{KS1}), marring an otherwise excellent study of
Leibnizian methodology \cite{Bos}, to which we now turn.

\section
{Grand programmatic statements}
\label{four}

In his seminal study, Bos argued that Leibniz exploited two competing
methods in his work, one \emph{Archimedean} and the other involving
the law of continuity and infinitesimals (see e.g.,
\cite[p.~57]{Bos}).
%
%here 57 should apparently be changed to 55
%

In asserting that Leibniz exploited distinct methods in developing the
calculus, we mean that he employed distinct conceptualisations of
continua; i.e., Leibniz employed different techniques for representing
relations among continuously changing magnitudes.  At a minimum, the
techniques differed in the inferences they sanctioned and in the
objects, whether ideal or fictional, which individual symbols in the
technique purported to represent.

Such a dichotomy can be reformulated in the terminology of dual
methodology (see \cite{KS1}) as follows.  One finds both A-track
(i.e., Archimedean) and B-track (Bernoullian, i.e., involving
infinitesimals)%
\footnote{ Scholars attribute the first systematic use of
infinitesimals as a foundational concept to Johann Bernoulli.  While
Leibniz exploited both infinitesimal methods and ``exhaustion''
methods (usually interpreted in the context of an Archimedean
continuum, but see footnote~\ref{exhaust}), Bernoulli never wavered
from the infinitesimal methodology.  To note the fact of such
systematic use by Bernoulli is not to say that Bernoulli's foundation
is adequate, or that that it could distinguish between manipulations
with infinitesimals that produce only true results and those
manipulations that can yield false results.}
methodologies in Leibniz.  In addition, Leibniz on occasion speculates
as to how one might seek to reformulate B-track techniques in an
A-track fashion.

Now there is no argument that such a pair distinct methodologies, A
and B, is present in Leibniz at the syntactic level.  Ishiguro does
not disagree with the apparent surface difference between them.  What
she argues, however, is that the syntactic difference is merely
skin-deep, so that once one clarifies the precise content of the
sentences one arrives at the conclusion that at that deeper level,
talk about infinitesimals (B-track) is merely shorthand for a
quantified statement (A-track), a position we denote
\[
B=A
\]
as shorthand for Ishiguro's contention that Leibnizian infinitesimals
are logical fictions. 

We argue that the syntactic difference in fact corresponds to a
semantic difference.  Each methodology has its respective ontology.
The B method involves a richer numerical structure than the A method.
Note that the structures have different ontological status.  The B
numerical structure involves pure fictions, while the A structure
involves ideal entities.  On this view, the A and B methods are truly
distinct, i.e., the Leibnizian infinitesimals are \emph{pure}
fictions, even though Leibniz occasionally argues that B should be
paraphrasable in terms of A, given enough effort.  This
\emph{hopefully paraphrasable} view can be denoted by the formula
\[
B>A,
\]
suggesting that what is involved in the B-method is an extended number
system including infinitesimals \emph{\`a la rigueur} (as Leibniz put
it with respect to ``des infinis'' in \cite[p.~92]{Le02a}), namely
what we refer to as a Bernoullian continuum.  D.~Jesseph expressed
this aspect of Leibniz's position in the following terms:
\begin{quote}
Leibniz often makes \emph{grand programmatic statements} to the effect
that derivations which presuppose infinitesimals can always be re-cast
as exhaustion proofs in the style of Archimedes. But Leibniz never, so
far as I know, attempted anything like a general proof of the
eliminability of the infinitesimal, or offered anything approaching a
universal scheme for re-writing the procedures of the calculus in
terms of exhaustion proofs. \cite[p.~233]{Je08}. [emphasis added]
\end{quote}
The most basic difference between the positions represented
respectively by~$B=A$ and~$B>A$ is that the former implies that the
background continuum of both the A-method and the B-method is
Archimedean, whereas the latter recognizes a genuinely enriched
(Ber\-noullian) continuum in the B-method.

We grant that for any given Leibnizian passage discussing the relation
between A-method and B-method, a plausible case can be made for
either~$B=A$ or~$B>A$, given sufficient ingenuity.  How can a scholar
determine which interpretation is truer to Leibniz's intentions?  In
the next few sections we present context-specific clues in Leibniz
that would allow one to choose between the two interpretations.

\section{Truncation manipulations}
\label{five}

Returning to the passage from the letter to Pinson quoted in
Section~\ref{three}, one discovers that Tho truncated the passage to
make it fit Ishiguro's analysis of Leibnizian infinitesimals.  The
full passage doesn't fit so well:
\begin{quote}
Car au lieu de l'infini ou de l'infiniment petit, on prend des
quantit\'es aussi grandes et aussi petites qu'il faut pour que
l'erreur soit moindre que l'erreur donn\'ee, de sorte que l'on ne
differe du style d'Archimede que dans les expressions qui sont plus
directes dans Nostre methode, et plus conforme \`a l'art
d'inventer.%
\footnote{We have retained Leibniz's spelling which differs slightly
from modern French spelling.  We provide an English translation as
found in \cite[note~25, p.~229]{Je08}: ``For in place of the infinite
or the infinitely small we can take quantities as great or as small as
is necessary in order that the error will be less than any given
error. In this way we only differ from the style of Archimedes in the
expressions, which are more direct in our method and better adapted to
the art of discovery.''}
\cite[p.~96]{Le01a}
\end{quote}
The conclusion of the passage, namely the clause concerning the
expressions
\begin{quote}
\emph{qui sont plus directes dans Nostre methode, et plus conforme \`a
l'art d'inventer}
\end{quote}
was omitted from Tho's translation cited in Section~\ref{three}.  This
conclusion clearly indicates that Leibniz's (B-track) method, where
the expressions are \emph{plus directes}, is distinct from the
(A-track) \emph{moindre que l'erreur donn\'ee} paraphrase thereof.
Leibniz's expression \emph{plus directes} suggests a distinct method
rather than merely a shorthand.  Thus Leibniz is following a distinct
strategy, that employs an enriched continuum.%
\footnote{\label{exhaust}In the context of Leibniz's reference to
Archimedes, it should be noted that there are other possible
interpretations of the exhaustion method of Archimedes.  The received
interpretation, developed in \cite{Di87}, is in terms of the limit
concept of real analysis.  However, in the 17th century, \cite
[pp.~280-290] {Wa85}
%Wallis
developed a different interpretation in terms of approximation by
infinite-sided polygons.  The ancient exhaustion method has two
components: (1) \emph{geometric construction}, consisting of
approximation by some \emph{simple} figure, e.g., a polygon or a line
build of segments, (2) \emph{justification} carried out in the theory
of proportion developed in \emph{Elements} Book V.  In the 17th
century, mathematicians adopted the first component, and developed
alternative justifications.}
The specific clues contained in this particular passage from Leibniz
favor the~$B>A$ reading over the~$B=A$ reading (see
Section~\ref{four}).

\section{Incomparables}

Let us examine Ishiguro's interpretation of Leibniz's notion of the
\emph{incomparably small}.

\subsection{Misleading and unfortunate}
\label{honor2}

Ishiguro claims that ``the incomparable magnitude is not an
infinitesimal magnitude" \cite[p.~87]{Is90}, and continues:
\begin{quote}
It is misleading for Leibniz to call these magnitudes incomparably
small.  What his explanation gives us is rather that a certain truth
about the existence of \emph{comparably} smaller magnitudes gives rise
to the notion of incomparable magnitudes, not incomparably smaller
magnitudes.  If magnitudes are incomparable, they can be neither
bigger nor smaller. \cite [p.~87--88]{Is90}
\end{quote}
She reiterates this claim in another sentence, regretting Leibniz's
choice of \emph{unfortunate} terminology:
\begin{quote}
As we have already mentioned, the unfortunate thing about Leibniz's
vocabulary here is that he moves from incomparable to incomparably
small or incomparably smaller (\emph{incomparabilitier parva} or
\emph{incomparabilitier minor}), when smaller is already a notion
involving comparison.  (ibid., p.~88)
\end{quote}
%{honor2}
For all her professed good intentions of defending Leibniz's honor
(see Section~\ref{honor1}), Ishiguro ends up being forced to defend
her interpretation by tarnishing that honor.  She reproaches him for
employing purportedly \emph{misleading} and \emph{unfortunate}
terminology in the context of incomparables.  In fact, Leibniz's
terminology for incomparables appears felicitous when the latter are
interpreted as \emph{pure} rather than logical fictions.

Ishiguro appears to claim that talk about \emph{incomparability}
excludes the relation of being smaller.  Note, however, that the term
\emph{incomparable} can be used in two distinct senses:
\begin{enumerate}
\item
it can refer to things that cannot be compared because they are of
different nature, e.g., a line and a surface;
\item
it can refer to things of the same nature but not comparable because
they are of different order of magnitude, e.g., an ordinary nonzero
real number and an infinitesimal.
\end{enumerate}
Ishiguro seems to assume the meaning (1) without any historical
evidence.  Seeking the meaning of incomparability she speculates
further:
\begin{quote}
The fact that we cannot add or subtract the quantities in question to
make one quantity constitutes, it seems, the very criterion of their
nonhomogeneity and hence of their incomparability.  (ibid., p.~88)
\end{quote}
Thus, on Ishiguro's reading,~$x$ is \emph{incomparably smaller}
than~$y$ means that~$x$ and~$y$ are not comparable at all, while their
\emph{incomparability} means that~$x$ cannot be added to~$y$.

\subsection{Evidence from \emph{Responsio}} 
\label{honor3}

To buttress her interpretation Ishiguro cites a few passages from
Leibniz.  The case of the 1695 \emph{Responsio} to Nieuwentijt is
particularly instructive here.  Ishiguro writes:
\begin{quote}
Leibniz, in the reply to Nieuwentijt of 1695 cited earlier, also
asserts that the magnitude of a line and a point of another line
cannot be added, nor can a line be added to a surface, and he says
that they are incomparable since only homogeneous quantities are
comparable. (Leibniz writes that all homogeneous quantities are
comparable in the Archimedean sense.) (ibid., p.~88)
\end{quote}
However, in the very same \emph{Responsio} we find a rather different
account of incomparability from Ishiguro's.  Leibniz writes as
follows:

\begin{quote}
[1] Caeterum aequalia esse puto, non tantum quorum differentia est
omnino nulla, sed et quorum differentia est incomparabiliter parva;
[2] et licet ea Nihil omnino dici non debeat, non tamen est quantitas
comparabilis cum ipsis, quorum est differentia.  [3] Quemadmodum si
lineae punctum alterius lineae addas, vel superficiei lineam,
quantitatem non auges.  [4] Idem est, si lineam quidem lineae addas,
sed incomparabiliter minorem.  [5]~Nec ulla constructione tale
augmentum exhiberi potest.  [6] Scilicet eas tantum homogeneas
quantitates comparabiles esse, cum Euclide lib.~5 defin.~5 censeo,
quarum una numero, sed finito%
\footnote{When Leibniz mentions a \emph{number} he is paraphrasing
Euclid's definition.  Now Euclid doesn't mention that he is speaking
of the line being added to itself \emph{finitely} many times.  But
this point is essential for Leibniz, as it isn't obvious that an
infinitesimal added to itself infinitely many times might not get to
be larger than a given finite magnitude.  Therefore Leibniz inserts
the term that Euclid doesn't use, and writes \emph{sed finito},
literally ``of course finite.''  Note that there is no \emph{number}
in Euclid's V.4, but rather a \emph{multitude}.  Thus, Leibniz reads
Euclid's \emph{multitude} as \emph{number}, and to be precise modifies
it by \emph{finite}.}
multiplicata, alteram superare potest.  [7] Et quae tali quantitate
non differunt, aequalia esse statuo, quod etiam Archimedes sumsit,
aliique post ipsum omnes.%
%\footnote{It is significant that in a book of essays called
%`Infinitesimal Differences. Controversies between Leibniz and his
%Contemporaries. Edited by Ursula Goldenbaum and Douglas Jesseph,
%(Walter de Gruyter, Berlin, New York 2007) there is no reference to
%this mention of Euclid's definition V.5 by Leibniz.}
%
\cite[p.~322]{Le95b} (numerals [1] through [7] added)
\end{quote}
%
%I think that those things are equal not only whose difference is absolutely nothing, but
%also whose difference is incomparably small; and although this difference need not be
%called absolutely nothing, neither is it a quantity comparable with those whose difference it
%is. Just as when you add a point of one line to another line or a line to a surface you do not
%increase the magnitude, it is the same thing if you add to a line a certain line, but one incomparably
%smaller. Nor can any increase be shown by any such construction. tr. by Douglas Jesseph . Strange enough, he quits at this point and do not porceede to the next sentence. In book of essays `Infinitesimal Differences. Controversies between Leibniz and his Contemporaries, Edited by Ursula Goldenbaum and Douglas Jesseph, Walter de Gruyter, Berlin, New York 2007 NO ONE refers to this mention of Euclid's definition V.5 by Leibniz
%
We translate this passage as follows:
\begin{enumerate}
\item[[1\!\!]] Furthermore I think that not only those things are
equal whose difference is absolutely zero, but also whose difference
is incomparably small.
\item[[2\!\!]] And although this [difference] need not absolutely be
called Nothing, neither is it a quantity comparable to those whose
difference it is.
\item[[3\!\!]] It is so when you add a point of a line to another line
or a line to a surface, then you do not increase the quantity.
\item[[4\!\!]] The same is when you add to a line a certain line that
is incomparably smaller.
\item[[5\!\!]] Such a construction entails no increase.
\item[[6\!\!]] Now I think, in accordance with Euclid Book V def.~5,
that only those homogeneous quantities one of which, being multiplied
by a finite number, can exceed the other, are comparable.
\item[[7\!\!]] And those that do not differ by such a quantity are
equal, which was accepted by Archimedes and his followers.
\end{enumerate}

Here Leibniz employs the term \emph{line} in the sense of what we
would today call a \emph{segment}.  In sentence [3], Leibniz exploits
the classical example with indivisibles (adding a point to a line
doesn't change its length) so as to motivate a similar phenomenon for
infinitesimals in sentence~[4] (adding an infinitesimal line to a
finite line does not increase its quantity), namely his law of
homogeneity (TLH) explained in more detail elsewhere (see
Section~\ref{three}).

In sentence [6], Leibniz refers to quantities satisfying Euclid~V.5,
i.e., the Archimedean axiom.  In a follow-up sentence [7], Leibniz
goes on to refer to `those [quantities],' say $Q$ and $Q'$, that `do
not differ by such a quantity,' namely they do \emph{not} differ by a
quantity of the type mentioned in sentence [6] (and \emph{satisfying}
Euclid~V.5 with respect to $Q$ or $Q'$).  Rather, $Q$ and $Q'$ differ
by a quantity \emph{not satisfying} Euclid~V.5.  Leibniz referred to
such quantities in sentence [1] as \emph{incomparably small}.  Thus
Leibniz is quite explicit about the fact that he is dealing with an
\emph{incomparably small} difference $Q-Q'$ which \emph{violates} the
Archimedean axiom V.5 when compared to either $Q$ or $Q'$.  Leibniz is
even more explicit about the fact that his \emph{incomparables}
violate Euclid V.5 (when compared to other quantities) in his letter
to l'Hospital from the same year:
\begin{quote}
J'appelle \emph{grandeurs incomparables} dont l'une multipli\'ee par
quelque nombre fini que ce soit, ne s\c cauroit exceder l'autre, de la
m\^eme facon qu'Euclide la pris dans sa cinquieme definition du
cinquieme livre. \cite[p.~288]{Le95a}
\end{quote}
The claim in \cite[p.~562]{Ar13} based on this very passage from
\emph{Responsio} that allegedly ``Leibniz was quite explicit about
this Archimedean foundation for his differentials as
`incomparables'\,'' is therefore surprising.  The 1695 letter to
l'Hospital does not appear in Arthur's bibliography.

Referring to the passage we quoted, Ishiguro claims that ``Leibniz, in
the reply to Nieuwentijt \ldots{} asserts that the magnitude of a line
and a point of another line cannot be added, nor can a line be added
to a surface, \ldots{}'' \cite[p.~88]{Is90}.  On the contrary, Leibniz
wrote in sentence~[3] that they can indeed be so added, though the
addition of a point to a line does not increase the line.  Thus
addition is possible according to Leibniz, contrary to what Ishiguro
claims Leibniz asserts.  It's just that according to Leibniz such
addition does not result in an increase.  What Leibniz actually wrote
undermines Ishiguro's claim about incomparables, rather than
supporting it.

Leibniz goes on to give a parallel example with infinitesimals in
sentence~[4].  Here addition is again possible, whereas its result is
unchanged in accordance with TLH (see Section~\ref{three}).  

Ishiguro somehow fails to mention the fact that Leibniz goes on to
give an example with infinitesimals.  This is not merely an instance
of truncation (see Section~\ref{five}).  Rather, it constitutes a
misrepresentation of Leibniz's position.

\subsection{Indivisibles, infinitesimals, and dimension}
\label{honor4}

Leibniz clearly understood the difference between infinitesimals (of
the same dimension as the quantities they modify) and indivisibles (of
positive codimension), contrary to Ishiguro's suggestion that
\begin{quote}
[t]he homogeneity of quantities in Leibniz \ldots{} seems not to
depend on a prior notion of a common dimension as in earlier
mathematicians, since Leibniz wanted to free mathematics from
geometrical intuitions.  \cite[p.~88]{Is90}
\end{quote}
The notion of ``common dimension'' is what distinguishes
infinitesimals from indivisibles.  Ishiguro's suggestion that Leibniz
did not distinguish between indivisibles and infinitesimals by means
of the notion of common dimension does no honor to Leibniz; see
Section~\ref{honor1}.

The passage in \cite[p. 322]{Le95b} cited above clearly indicates what
Leibniz means by comparable quantities.  Namely,~$x$ and~$y$ are
comparable when the following condition is satisfied:
\[
(\exists n\in \N)(nx>y),\] i.e.~$x, y$ do obey Euclid's definition
V.4%
\footnote{\label{V4}Leibniz lists number V.5 for Euclid's definition
instead of V.4.  In some editions of the \emph{Elements} this
definition does appear as V.5.  Thus, \cite{Eu} as translated by
Barrow in 1660 provides the following definition in V.V (the notation
``V.V'' is from Barrow's translation): \emph{Those numbers are said to
have a ratio betwixt them, which being multiplied may exceed one the
other}.}
as cited in Leibniz's sentence~[6].  Leibniz defines an {incomparably
small} in terms of a violation of~V.4.  Thus, even though Leibniz
eschews geometrical intuition, he is still able to distinguish
indivisibles from infinitesimals.

\subsection{Theory of magnitudes}

Since Leibniz explicitly refers to Euclid's definition V.4 in the
\emph{Responsio} (see Subsection~\ref{honor3}, sentence~[6]), let us
turn to the theory of magnitudes as developed in Book V of the
\emph{Elements}.  Euclid's magnitudes of the same kind (homogeneous
quantities in Leibniz's terminology) can be formalized as an ordered
additive semigroup with a total order,~$\mathfrak M=(M,+,<)$,
characterized by the five axioms given below.  

\cite{Bec} and \cite[pp.~101-122]{BM} provide detailed sources for the
axioms below in the primary source (Euclid).  See also
\cite[pp.~118-148]{Mue} which mostly follows Beckmann's development.
Axiom E1 below interprets Euclid~V.4:

\begin{enumerate}
\item[E1] $(\forall x, y\in M)(\exists n\in \N)(nx > y)$,
\item[E2] $(\forall x, y\in M)(\exists z\in M)(x < y \Rightarrow x + z=y)$,
\item[E3] $(\forall x, y, z\in M)(x < y \Rightarrow x + z < y + z)$,
\item[E4] $(\forall x\in M)(\forall n \in \N)(\exists y\in M)(x =ny)$,
\item[E5] $(\forall x, y, z\in M)(\exists v\in M)(x : y :: z : v).$
\end{enumerate}
Comparable quantities can both be added to one another, and they are
also subject to the relations \emph{greater than} and \emph{less
than}.  It follows from these axioms that for any~$x,y\in M$ the
following inequality holds:
\[
y<y+x.
\]
In the realm of incomparable quantities this inequality does not hold,
even though incomparables can be added.  Leibniz's claim can be
formalized as the relation $y\not< y+x$ characterizing incomparable
quantities.%

\subsection{\emph{Elements} Book VI on horn angles}

We turn next to Ishiguro's claim that incomparable quantities cannot
be compared at all by means of the relations \emph{greater than} or
\emph{less than}.

In Book VI of the \emph{Elements} one finds that line segments form a
semigroup of magnitudes of the same kind,~$\mathfrak M_1$, triangles
form another,~$\mathfrak M_2$, rectilinear angles form yet
another,~$\mathfrak M_3$;%
\footnote{See e.g., \cite{Euc}, VI.1, 2, 33.}
there are other kinds of magnitudes in addition to the ones just
mentioned.  Euclid deals with two kinds of angles in the
\emph{Elements}: the first kind consists of rectilinear angles, while
the second kind consists of angles cut out/formed by a line and an arc
of a circle.%
\footnote{In the 17th century such angles were called \emph{horn
angles}.}
These two kinds of angles are compared in proposition~III.16. Its
thesis reads as follows:

\begin{quote}
A (straight-line) drawn at right-angles to the diameter of a circle,
from its end, will fall outside the circle. And another straight-line
cannot be inserted into the space between the (aforementioned)
straight-line and the circumference. And the angle of the semi-circle
is greater than any acute rectilinear angle whatsoever, and the
remaining (angle is) less (than any acute rectilinear angle).
(translation by Fitzpatrick in \cite{Euc})
\end{quote}
We present below the accompanying diagram.
%
%the source pdf is in texfiles/horowitz/pdfbackup/Tw16_ks3.pdf

\centerline{\includegraphics[scale=0.5]{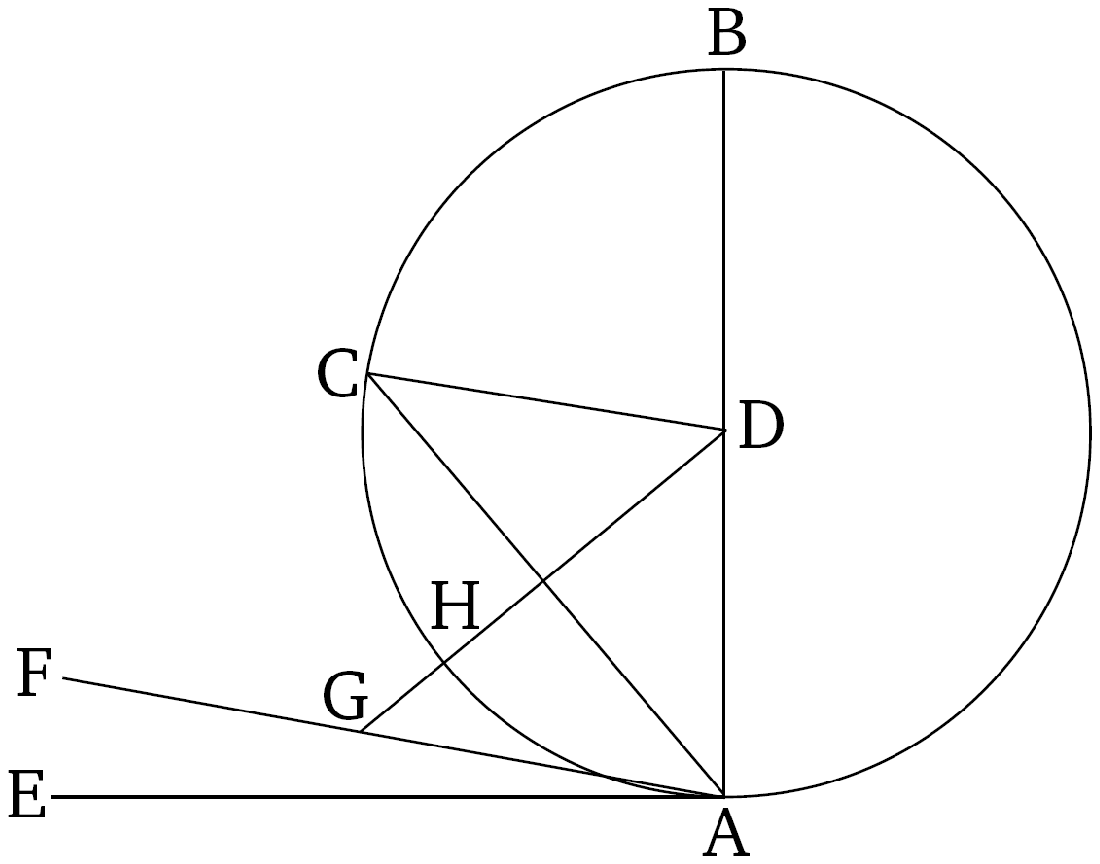}} 

Here ``the remaining" angle, that is the one formed by the arc~$CA$
and the tangent line~$EA$, does not belong to the kind (i.e., species)
of rectilinear angles~$\mathfrak M_3$.  Euclid proves it to be less
than any acute rectilinear angle.%
\footnote{It is worth noting that there are angles cut out by a line
and a curve in Leibniz papers and the phrase `infinitely small angle'
occurs many times; see for example \cite{Ch}.}

From the point of view of Greek mathematics, one can construct
incomparable quantities~$x, y$, meaning that they are not of the same
kind and do not obey the Archimedean axiom, while at the same time the
relation~$x<y$ obtains.  Here~``$x$ is incomparably smaller than~$y$''
means~$x$ is smaller than~$y$ and~$x, y$ are incomparable, which can
be formalized as follows:
\[
x<y\ \ \mbox{and}\ \ y\not< y+x.
\]
Thus, incomparable quantities can be compared by inequalities both
according to Euclid and according to Leibniz.

\section{Textual evidence}
\label{text}

In this section we will examine the textual evidence Ishiguro presents
to support her claim that infinitesimals are logical fictions.

\subsection{Letter to Varignon}
\label{71}

Ishiguro's first piece of textual evidence in favor of her logical
fiction hypothesis is a letter to Varignon dated~2~february 1702
\cite{Le02a}.  Ishiguro (p.~82) does not provide a direct quotation
but refers to \cite{Ge50}, vol.~IV, p.~93, which contains what seems
to be one of two occurrences in Leibniz of the term ``syncategorematic
infinite''.%
\footnote{The other one is in Leibniz's correspondence with des
Bosses: Gerhardt (ed.), Leibniz, Philosophische Schriften vol.~II,
pp.~314f.: Datur infinitum syncategorematicum etc.}

The 2~february~1702 letter exploits the term \emph{syncategorematic}.
However, it is not obvious that Leibniz uses it in the same technical
sense as Ishiguro.  Leibniz discusses a number of examples including
imaginary numbers, dimensions beyond~$3$, and exponents which are not
\emph{ordinary numbers}, and then comments as follows:
\begin{quote}
Cependant il ne faut pas s'imaginer que la science de l'infini est
degrad\'ee par cette explication et reduite \`a des fictions; car il
reste tousjours un infini \emph{syncategorematique}, comme parle
l'ecole, \dots{} \cite[p.~93]{Le02a}. [emphasis added]
\end{quote}
Leibniz then goes on to discuss the summation of a geometric series
and points out that no infinitesimals need appear here.  He is
discussing a way of accounting for B-methodology in terms of
A-methodology.  The plain meaning of the text, as mentioned in
Section~\ref{one}, is that there is a pair of distinct methodologies,
and if the \emph{fictions} of the B method were found lacking, one
could, at least in principle (recall Jesseph's remark concerning
\emph{grand programmatic statements}), fall back on an A-type
\emph{syncategorematic} paraphrase.

In analyzing this occurrence of the adjective \emph{syncategorematic},
we again have the problem of investigating which of the two
interpretations,~$B=A$ or~$B>A$, is more faithful to Leibniz's general
philosophical outlook (see Section~\ref{four}).  We will therefore
look for additional clues in the letter that may favor one of the
interpretations.

It is significant that the letter also contains a discussion of the
\emph{law of continuity}.  Here Leibniz writes that the rules of the
finite succeed in the infinite, and vice versa:
\begin{quote}
\ldots{} il se trouve que \emph{les r\`egles du fini r\'eussissent
dans l'infini} comme s'il y avait des atomes (c'est \`a dire des
\'el\'ements assignables de la nature) quoiqu'il n'y en ait point la
mati\`ere \'etant actuellement sousdivis\'ee sans fin; \emph{et que
vice versa les r\`egles de l'infini r\'eussissent dans le fini}, comme
s'il y'avait des infiniment petits m\'etaphy\-siques, quoiqu'on n'en
n'ait point besoin; \cite[p.~93-94]{Le02a} [emphasis added]
\end{quote}
Leibniz goes on to mention the \emph{souverain principe}:
\begin{quote}
et que la division de la mati\`ere ne parvienne jamais \`a des
parcelles infiniment petites: c'est parce que tout se gouverne par
raison, et qu'autrement il n'aurait point de science ni r\`egle, ce
qui ne serait point conforme avec la nature du \emph{souverain
principe} \cite[p.~94]{Le02a}. [emphasis added]
\end{quote}
A number of scholars, including \cite[p.~67]{Kn02} as well as
\cite[p.~145]{Lau92}, identify the passage on p.~93-94 as an
alternative formulation of the law of continuity, viz., \emph{the
rules of the finite succeed in the infinite, and conversely}.  Thus,
recent scholarship has interpreted this passage as Leibniz's
endorsement of the possibility of transfering properties from finite
numbers to infinite (and infinitesimal) numbers and vice versa.  For
example, the usual rules governing the arithmetic operations and
elementary functions should be obeyed by infinitesimals, as well.

Now if infinitesimal expressions were merely shorthand for talk about
ordinary finite numbers or a sequence thereof, Leibniz's law of
continuity would amount to the assertion that
\begin{quote}
\emph{each element in a sequence of ordinary numbers obeys the same
rules as ordinary numbers}.
\end{quote}
But this seems anticlimactic, and moreover too tautological to have
been termed a law or a \emph{souverain principe} by Leibniz.  Leibniz
writes further:
\begin{quote}
Et c'est pour cet effect que j'ay donn\'e un jour des lemmes des
incomparables dans les Actes de Leipzic, qu'on peut entendre comme on
vent [sic], soit des infinis \emph{\`a la rigueur}, soit des grandeurs
seulement, qui n'entrent point en ligne de compte les unes au prix des
autres. \cite[p.~92]{Le02a}. [emphasis added]
\end{quote}
Leibniz's pair of ``soit"s in this remark indicates that there is a
pair of distict methodologies involved, as we elaborated in
Section~\ref{four}.  Note that Ishiguro quotes both a passage on page
91 preceding this remark \cite[p.~86]{Is90}, and a passage on page 92
following the remark \cite[p.~87]{Is90}, but fails to quote this
crucial remark itself (see our Section~\ref{five} on truncation
manipulations).

Thus, the letter offers support for the thesis that Leibniz thought
infinitesimals (and infinite numbers) could stand on their own
(\emph{\`a la rigueur}), without paraphrase in terms of finite
quantities.  The letter fits well with the~$B>A$ idea that on occasion
Leibniz tried to argue optimistically (see Section~\ref{four} on
\emph{grand programmatic statements}) that B-track techniques should
be paraphrasable in terms of A-track ones.

\subsection{\emph{Theodicy}}

Ishiguro's second piece of textual evidence is from Leibniz's
\emph{Theodicy}:
\begin{quote}
``every number is finite and assignable, every line is also finite and
assignable.  Infinites and infinitely small only signify magnitudes
which one can take as big or as small as one wishes, in order to show
that the error is smaller than the one that has been assigned''
(\emph{Theodicy}, \S 70).  \cite[p.~83]{Is90}.
\end{quote} 
We have retained Ishiguro's precise punctuation including the
quotation marks.  Note that her quotation marks close the citation
without any punctuation mark at the end of the citation.  There is no
indication in Ishiguro that Leibniz's sentence does not end there, but
rather continues.  It is instructive to examine Leibniz's sentence in
full:
\begin{quote}
Every number is finite and specific; every line is so likewise, and
the infinite or infinitely small signify only magnitudes that one may
take as great or as small as one wishes, to show that an error is
smaller than that which has been specified, that is to say, that there
is no error; or else by the infinitely small is meant the state of a
magnitude at its vanishing point or its beginning, conceived after the
pattern of magnitudes already actualized. \cite{Le10a} (translation by
Gutenberg Project)
\end{quote}
The closing phrase,
\begin{quote}
\emph{or else by the infinitely small is meant the state of a
magnitude at its vanishing point or its beginning, conceived after the
pattern of magnitudes already actualized,}
\end{quote}
was truncated by Ishiguro.  In this omission she is not without
co-conspirators: the same \emph{truncated} passage appears in
\cite[p.~76]{Go08} and \cite[p.~29]{Fe08}.
%
%Ferraro cites the Dutens edition 1:107; see commented out entry in
%bibliography
%

Leibniz's conclusion in \S\,70 suggests that there does exist a way of
working with infinitesimals \emph{\`a la rigueur}.  This would
presumably involve an enriched system of magnitudes, whose additional
elements shared properties with the (\emph{already actualized})
elements in the original system.  Leibniz is being rather vague here
and it is hard to know what he means exactly by magnitudes being
``conceived after the pattern of magnitudes already actualized,''
especially since \S\,70 is preceded by \S\,69 on free will and
followed by~\S\,71 on the Gospels, making it difficult to rely on the
context for a clarification.  However, our impression is that Ishiguro
is not telling the full story here, for she observes:
\begin{quote}
As the \emph{Theodicy} is a very late book (1710), it may be
\emph{thought} that this expresses a later-year shift to finitism
brought about by \emph{senility}.  In order to see that this is not
the case, let us trace some of the things Leibniz wrote on
infinitesimals from his early years.  \cite [p.~83] {Is90} [emphasis
added]
\end{quote}
To be sure, Ishiguro soon enough rejects her \emph{senility}
hypothesis.  However, even a hypothesis that is ultimately rejected
must have a grain of plausibility to it.  Otherwise why would one want
to consider it in the first place?  Note that the Leibniz-Clarke
correspondence is well-regarded, and it comes at the very end of
Leibniz's life.  Leibniz died in 1716 before he had a chance to
respond to Clarke's 5th letter.

Ishiguro's \emph{thought} here seems at odds with her stated goal of
defending Leibniz's honor (see Section~\ref{one}).  Her \emph{thought}
seems to imply that Leibniz shifted to an infinitesimal-barring
finitism in 1710.  The truth is that, on the contrary, Leibniz was at
that time at the height of his intellectual powers, and was as
committed as ever to developing the B-methodology, including its
foundations, as is evidenced by his extremely lucid 1710 text on the
transcendental law of homogeneity \cite{Le10b} analyzed in
Section~\ref{three}.

\subsection{\emph{Nova Methodus}}
\label{six}

Ishiguro's third piece of textual evidence is drawn from \emph{Nova
Methodus} \cite{Le84} (a text she misdates at 1685).  She makes
several dubious claims related to this work.

First, she alleges that in this text, differentials are ``defined
through the proportion of finite line segments'' \cite[p.~83]{Is90}.
What Leibniz actually writes is the following:
\begin{quote}
Now some right line taken arbitrarily may be called~$dx$, and the
right line which shall be to~$dx$, as~$v$ (or~$w$,~$y$,~$z$, resp.)
is to~$VB$ (or~$WC$,~$YD$,~$ZE$, respect.) may be called~$dv$ (or
$dw$,~$dy$,~$dz$, resp.), or the differentials \ldots{}
\cite[p.~467]{Le84}.
\end{quote}
This passage in a notoriously (and deliberately) obscure work cannot
qualify as a definition of \emph{differential}, and certainly offers
no support for Ishiguro's claim that infinitesimal expressions are
non-referring.  Leibniz scholars have argued that he had to conceal
the use of infinitesimals in this publication to avoid the wrath of
opponents:
\begin{quote}
The structure of the text [i.e., \emph{Nova Methodus}], which was much
more concise and complex than the primitive Parisian manuscript
essays, was complicated by the \emph{need to conceal the use of
infinitesimals}.  Leibniz was well aware of the possible objections he
would receive from mathematicians linked to classic tradition who
would have stated that the infinitely small quantities were not
rigorously defined, that there was not yet a theory capable of proving
their existence and their operations, and hence they were not quite
acceptable in mathematics. \cite[p.~49]{Ro05} [emphasis added]
\end{quote}
This would account for the obscurities of Leibniz's discussion of
differentials here, which offers no support at all for a
syncategorematic reading of Leibnizian infinitesimals.

Furthoremore, Ishiguro goes on to provide a syncategorematic
interpretation of Leibniz's construction of a line through a pair of
infinitely close points:
\begin{quote}
Leibniz writes that a tangent is found to be a straight line drawn
between two points on a curve of infinitely small distance, or a side
of a polygon of infinite angles.  However, \ldots{} infinitely small
distances can be thought of as distances that can be taken smaller
than any distances that are given. \cite[p.~84]{Is90}.
\end{quote}
This passage furnishes no explanation for the asymmetry of the two
points involved in the received definition of the tangent line via
secant lines (as discussed in Section~\ref{eight}).  

As we already mentioned in Section~\ref{three}, \emph{Nova Methodus}
contains the first mention of Leibniz's \emph{law of homogeneity},
evidence in favor of infinitesimals \emph{\`a la rigueur}.  Thus, the
clues contained in \emph{Nova Methodus} support the~\hbox{$B>A$}
reading.

\subsection{\emph{Responsio a Nieuwentijt}}

Ishiguro's fourth piece of textual evidence is the 1695 response to
Nieuwentijt \cite{Le95b} published in \emph{Acta Eruditorum}.  She
writes that
\begin{quote}
Leibniz explains that although he treats (\emph{assumo}) infinitely
small lines~$dx$ and~$dy$ as true quantities sui generis, this is just
because he found them useful for reasoning and discovery.  I take it
that he is treating them as convenient \emph{theoretical fictions}
because using signs which looks [sic] as if they stand for quantities
sui generis is useful.  \cite[p.~84]{Is90} [emphasis added]
\end{quote}
In point of fact, \emph{theoretical fictions} are on a par with what
we refer to as \emph{pure fictions}.  What Ishiguro writes here
undermines her syncategorematic interpretation and supports ours.  The
fact that this is what the passage means is demonstrated by the
comparison with imaginaries, for which Leibniz has no syncategorematic
account.  The passage Ishiguro is referring to reads as follows:
\begin{quote}
Itaque non tantum lineas infinite parvas, ut~$dx$,~$dy$, pro
quantitatibus veris in suo genere assumo, sed et earum quadrata vel
rectangula~$dxdx$,~$dydy$,~$dxdy$, idemque de cubis aliisque
altioribus sentio, praesertim cum eas ad ratiocinandum inveniendumque
utiles reperiam.%
\footnote{This is translated as follows by Parmentier: ``Ainsi au
nombre des grandeurs r\'eelles en leur genre, je ne compte pas
seulement les lignes infiniment petites~$dx$,~$dy$, mais aussi leurs
carr\'es ou leurs produits~$dxdx$,~$dxdy$,~$dydy$, il en va de m\^eme
d'apr\`es moi de leurs cubes et de leurs puissances sup\'erieures,
compte tenu notamment de la f\'econdit\'e que j'y ai d\'ecouverte dans
les raisonnements et les inventions.'' \cite[p.~328]{Le1989}.}
\hfil\linebreak \cite[p.~322]{Le95b}.
\end{quote}
According to Ishiguro, Leibniz says that he treats~$dx$ and~$dy$ as
true quantities sui generis \emph{just because} he found them useful.
But is that really what he is saying?  Leibniz refers to
infinitesimals by the adjective \emph{veris} meaning ``true'' or
``real''.  If his infinitesimals were logical fictions i.e., merely
shorthand for sequences of real values, what novelty would there be in
emphasizing, as he does, that he includes infinitesimals (as well as
those of higher order) among what he describes as true or real
quantities?  Why emphasize this point if infinitesimals were merely
shorthand for sequences of what are already ordinary values drawn from
an Archimedean system?  Furthermore, why would he seek to buttress
such a straightforward point by underscoring the usefulness of
infinitesimals in reasoning and discovery?

A few lines earlier on page 322, Leibniz cites Euclid V.5%
\footnote{This corresponds to V.4 in modern editions; see
footnote~\ref{V4}.}
in a way similar to the 1695 letter to l'H\^opital, indicating a
violation of the Archimedean property; see Section~\ref{eight}.
Remarkably, Leibniz uses the term \emph{numerus infinitus}, meaning
infinite \emph{number}--rather than infinite quantity--here, blocking
the option of interpreting it as a variable \emph{quantity} increasing
without bound.%
\footnote{Parmentier in his 1989 French translation devotes a lengthy
footnote 30 on page~325 to Leibniz's usage of the term \emph{numerus}
here.}

Leibniz not only speaks of two distinct methods but gives them names
that suggest what his personal preferences are.  Namely, Leibniz
describes what we refer to as the A-method as \emph{reducendi via}
(the way of reducing)%
\footnote{This was rendered \emph{cette d\'emonstration r\'egressive}
in Parmentier's translation \cite[p.~327]{Le1989}.}
and the infinitesimal method as \emph{methodus directa} (the direct
method).  It is instructive to analyze the relevant passage in detail.
Leibniz writes in this \emph{Responsio}:
\begin{quote}
Quoniam tamen methodus directa brevior est ad intelligendum et utilior
ad inveniendum, sufficit cognita semel reducendi via postea methodum
adhiberi, in qua incomparabiliter minora negliguntur, quae sane et
ipsa secum fert demonstrationem suam secundum lemmata a me Febr. 1689
communicata. \cite{Ge50}, vol.~V, p.~~322.
\end{quote}
We translate this passage as follows:
\begin{quote}
But since the direct method [\emph{methodus directa}] is shorter to
understand and a more useful way of finding [i.e., discovering], it
suffices, once the way of reducing [\emph{reducendi via}] is known, to
apply afterward the method in which quantities that are incomparably
smaller are neglected, which in fact carries its own demonstration
according to the lemmas that I communicated in February, 1689.
\end{quote}
What emerges from this sentence is that there are two distinct
methods, [A] ``via reduction" and [B] a ``direct method" using
infinitesimals.  The infinitesimal method is riskier but more
powerful, and what Leibniz is pointing out is that having gained some
experience with the traditional method so that one already knows what
kind of results to expect, one can safely used the infinitesimal
method that yields the same results but more efficiently.  Leibniz
points out that once the reductive method A has been used and
understood, from that point onward one can systematically use the
direct method B (which involves discarding infinitesimals), since it
is quicker and more useful.  These clues furnish further evidence in
favor of the~$B>A$ reading over the~$B=A$ reading (see
Section~\ref{four}).

\subsection{The 7 june 1698 letter to Bernoulli}
\label{7june}

Ishiguro's fifth piece of textual evidence is the 7 june 1698 letter
to Bernoulli \cite{Le98}.%
\footnote{Ishiguro gives the page range as 499-500.  Actually the
letter occupies four pages 497-500.}
She writes that Leibniz likens the status of infinitesimals to that of
imaginary numbers in this letter (ibid., p.~84).

Since Ishiguro does not elaborate any further, it is difficult to see
how this could be interpreted as a piece of evidence in favor of her
logical fiction hypothesis, since in point of fact complex numbers
could not (in Leibniz's day) be replaced by quantified paraphrases
ranging over ordinary numbers, so complex numbers (or imaginary
quantities, as Leibniz called them) are \emph{pure fictions} par
excellence.  Leibniz repeatedly insisted (not merely in this letter to
Bernoulli) on the analogy between the fictional status of
infinitesimals and complex numbers.  Meanwhile, Leibniz described
imaginaries as having their \emph{fundamentum in re} (basis in fact)
in \cite[p.~93]{Le95a}.  The comparison to complex numbers tends to
undermine the \emph{logical fiction} hypothesis concerning Leibnizian
infinitesimals.  This theme was explored more fully in \cite{SK}.

\section{Euclid V.4, Apollonius, and tangent line}
\label{eight}

According to the letter to l'Hospital \cite{Le95a}, Leibniz's
infinitesimals violate Euclid~V.4:%
\footnote{Leibniz actually refers to V.5; see footnote~\ref{V4}.}
\begin{quote}
J'appelle \emph{grandeurs incomparables} dont l'une multipli\'ee par
quelque nombre fini que ce soit, ne s\c{c}auroit exceder l'autre, de
la m\^eme fa\c{c}on qu'Euclide la pris dans sa cinquieme definition du
cinquieme livre. \cite[p.~288]{Le95a}
\end{quote}
Note Leibniz's use of the term \emph{grandeur}, i.e., magnitude,
rather than the more ambiguous term \emph{quantity}.  A magnitude
(e.g., 5 feet) is a level of a quantity (length).  Here the option of
interpreting this as shorthand for \emph{variable} quantity is not
available, barring also a \emph{logical fiction} reading.  The
definition Leibniz refers to is a variant of what is known today as
the Archimedean property of continua.  This indicates that Leibniz
embraces what we refer to as a Bernoullian continuum (though certainly
not a non-Archimedean continuum in a modern set-theoretic sense),
contrary to Ishiguro's Chapter~5.

%Jesseph
%
\cite{Je15} shows that strategies Leibniz employed in the attempt to
show that such fictions are acceptable because the use of
infinitesimals can ultimately be eliminated have to presume the
correctness of an infinitesimal inference (i.e., inference exploiting
infinitesimals), namely identifying the tangent line to a curve as
part of the construction.  In the case of conic sections this strategy
succeeds because the tangents are already known from Apollonius.  But
for general curves (including transcendental ones treated by Leibniz),
infinitesimals \emph{\`a la rigueur} remain an irreducible part of the
Leibnizian framework, contrary to Ishiguro's Chapter~5.

In 1684, Leibniz wrote as follows concerning the tangent line:
\begin{quote}
\ldots{} to find a tangent is to draw a straight line, which joins two
points of the curve having an infinitely small difference \ldots{}
\cite{Le84}
\end{quote}
The definition of a tangent line as the line through a pair of
infinitely close points on the curve poses a challenge to a
proto-Weier\-strassian reading.  Such a reading involves having to fix
one of the points and to vary the other, and construct a sequence of
secant lines producing the tangent line in the limit.  In such a
reading, one of Leibniz's points would be a genuine mathematical
concept (the future point of tangency), while the other, merely a
syncategorematic device or a shorthand for a \emph{sequence} of
ordinary values.  

However, nothing whatsoever about Leibniz's wording would indicate
that there is such an asymmetry between the two points, and on the
contrary implies a symmetry between them: either both denote, or
neither denotes.  Leibniz's definition of the tangent line is at odds
with Ishiguro's Chapter~5.

The most devastating blow to Ishiguro's Chapter~5 is the
\emph{hierarchical} structure on the
Leibnizian~$dx$'s,~$dx^2$'s,~$ddx$'s, etc., ubiquitous in Leibniz's
texts.  One can replace~$dx$ by a sequence of finite
values~$\epsilon_n$ and furnish a concealed quantifier incorporated
into a hidden proto-Weierstrassian limit notion so as to
interpret~$dx$ as shorthand for the sequence~$\langle \epsilon_n :
n\in\mathbb{N}\rangle$.  However, one notices
that~$\lim_{n\to\infty}\epsilon_n=0$, as well
as~$\lim_{n\to\infty}\epsilon_n^2=0$, and also
unsurprisingly~$\lim_{n\to\infty}(\epsilon_n+\epsilon_n^2)=0$.  Thus,
the Leibnizian substitution~$dx+dx^2=dx$ in accordance with the TLH
becomes a meaningless tautology~$0+0=0$.  To interpret it in a
meaningful fashion, Ishiguro would have to introduce additional ad hoc
proto-Weierstrassian devices%
\footnote{E.g., for every positive epsilon there exists a positive
delta such that whenever~$dx$ is less than delta, the difference
$\left|\frac{dx+dx^2}{dx} - 1\right|$ is less than epsilon.}
with no shadow of a hint in the original Leibniz.

\section{Conclusion}

Leibniz on occasion writes that arguments using infinitesimals
(B-track terminology) could be paraphrased in terms of ordinary
numbers drawn from an Archimedean number system (A-track terminology).
The question we have investigated is what exactly is involved in such
a paraphrase.  Ishiguro argued that Leibnizian infinitesimals do not
designate, so that when one clarifies the logical content of his
propositions mentioning infinitesimals, the infinitesimals disappear
and one is left with a suitable quantified proposition.  Ishiguro's
claim is that Leibnizian infinitesimals are \emph{logical fictions}.
We have argued that Leibnizian infinitesimals are \emph{pure fictions}
not eliminable by paraphrase.

%Leibniz's \emph{grand programmatic statements} (see
%Section~\ref{four}) concerning infinitesimal arguments being amenable
%to paraphrase, can be clarified in light of a modern metatheorem
%asserting that every result provable by hyperreal techniques, also
%admits a proof purely in the real framework (though such a proof may
%be vastly more complicated).\footnote{\label{nelson}Nelson's
%conservation theorem \cite[Section 8] {Ne} asserts that if an ordinary
%mathematical statement about the reals is proved in IST then it can be
%also proved in ZFC.  Here ZFC is the Zermelo--Fraenkel axiomatic set
%theory commonly accepted as the axiomatic base of modern mathematics,
%and IST is most commonly used extension of ZFC by \emph{nonstandard}
%methods which include arguments with infinitesimals and infinitely
%large numbers.  See \cite{KR} for a detailed study of the
%conservativity and other metamathematical properties of nonstandard
%set theories.  However if one restricts the accepted toolkit of
%mathematical methods to a theory~$T$ substantially weaker than ZFC
%then it may happen that the corresponding nonstandard extension~$T^*$
%of~$T$ is strictly stronger than~$T$ itself, as demonstrated in
%\cite{HK}.  This applies, for instance, in the case when~$T$ is
%the~$n$th order Peano arithmetic PA.  A corresponding nonstandard
%theory~$T^*$ turns out to be as strong as the~$(n+1)$th order PA,
%hence, stronger than~$T$ itself.}  

This does not mean that Leibniz's infinitesimals are Robinson's
infinitesimals; far from it.  The well-known differences between them
(Leibniz's continuum being nonpunctiform whereas Robinson's,
punctiform) should be approached from the viewpoint of the distinction
between \emph{mathematical practice} and \emph{the ontology of
mathematical entities} developed in \cite{Be65}, \cite{Qu}.  What
emerges from our analysis is that modern infinitesimal frameworks
provide better proxies for understanding Leibnizian procedures and
actual mathematical \emph{practice} than the Weierstrassian framework
(similarly punctiform, like Robinson's) Ishiguro seeks to read into
Leibniz.

Ishiguro's syncategorematic reading is contrary to explicit Leibnizian
texts, such as his 1695 texts \emph{Responsio} to Nieuwentijt and
letter to l'Hopital where he writes that his differentials violate
Euclid Def.~V.4, closely related to the Archimedean property of
continua.  Leibniz describes B-track methods as being \emph{direct}
and A-track methods as involving (indirect) \emph{reductio} arguments,
implying distinct methodologies.  Leibniz repeatedly likens
infinitesimals to imaginaries, and at least once described the latter
as having their \emph{fundamentum in re} (basis in fact), suggesting
that both are entities.  In some cases, Ishiguro resorts to
misrepresentation of what Leibniz wrote so as to buttress her position
(see Subsection~\ref{honor3} on the possibility of addition of
incomparables).

In view of all the difficulties with Ishiguro's reading, we can only
conclude that the legitimate grounds for a ``rehabilitation'' (if any
is needed) of Leibniz's infinitesimal calculus are to be found in the
Leibnizian theory itself%
\footnote{``If the Leibnizian calculus needs a rehabilitation because
of too severe treatment by historians in the past half century, as
Robinson suggests (1966, p. 250), I feel that the legitimate grounds
for such a rehabilitation are to be found in the Leibnizian theory
itself.''  \cite[pp.~82-83]{Bos}}
(including his transcendental law of homogeneity), rather than in
either Fregean quantifiers, Weierstrassian epsilontics, or Russellian
logical fictions.

\section*{Acknowledgments}

V. Kanovei received partial support from the RFBR grant 13-01-00006.
M.~Katz was partially funded by the Israel Science Foundation grant
number~1517/12 and RSF grant 14-50-00150.  We are grateful to Thomas
Mormann and to the anonymous referees for helpful suggestions.  The
influence of Hilton Kramer (1928-2012) is obvious.

\end{document}